\input amstex
\documentstyle{amsppt}


\hcorrection{0.75truein}

\def\C{{\Cal C}}

\def\HH{{\Bbb H}}
\def\RR{{\Bbb R}}

\def\card{\operatorname{card}}
\def\kbar{\overline{k}}

\def\ktilde{\widetilde{k}}
\def\Rtilde{\widetilde{R}}

\def\ATV{{\bf ATV}}
\def\Co{{\bf Co}}
\def\CH{{\bf CH}}
\def\EisRo{{\bf ER}}
\def\FS{{\bf FS}}
\def\GS{{\bf GS}}
\def\Jat{{\bf Ja}}
\def\Lam{{\bf La}}
\def\Len{{\bf Le}}

\def\McPet{{\bf MP}}
\def\McRob{{\bf MR}}
\def\Pas{{\bf Pa}}
\def\Pie{{\bf Pi}}
\def\Rob{{\bf Rb}}
\def\Rwn{{\bf Rw}}
\def\ssw{{{\bf SW1}}}
\def\sww{{\bf SW2}}

\def\Yam{{\bf Ya}}
\def\Va{{\bf Va}}
\topmatter

\title Simplicity of noncommutative Dedekind domains \endtitle

\author K. R. Goodearl and J. T. Stafford \endauthor

\address Department of Mathematics, University of California, Santa
Barbara, Ca 93106-3080\endaddress
\email goodearl\@math.ucsb.edu \endemail
\address Department of Mathematics, University of Michigan, Ann Arbor,
MI 48109-1109 \endaddress
 \email jts\@umich.edu \endemail

\thanks The research of both authors was partially supported by grants
from the National Science Foundation. Some of it was carried out while
the authors participated in the Noncommutative Algebra Year (1999-2000)
at the Mathematical Sciences Research Institute in Berkeley, and they
thank MSRI for its support.
\endthanks

\abstract The following dichotomy is
established: A finitely generated, complex
  Dedekind domain that is {\it not} commutative is a simple
ring. Weaker versions of this dichotomy are proved for Dedekind prime rings and
hereditary noetherian prime rings.
\endabstract



\endtopmatter

\document

\head Introduction \endhead
When the classical concept of a Dedekind domain was extended to
noncommutative rings, the natural examples that arose were either
classical orders (and hence finitely generated modules over their centres)
 or simple rings such as the Weyl algebra
$A_1({\Bbb C})$. Indeed, among finitely generated Dedekind domains over
algebraically closed fields,
 classical orders and simple rings are the only known examples.
 This dichotomy in the examples
suggests that an actual dichotomy might exist among general Dedekind domains,
although we are not aware that any such conjecture has been formulated in the
literature. The main goal in this paper is to establish just such a result
as well as give similar dichotomies for  Dedekind prime rings and HNP rings.

Before stating the first theorem we need some definitions. An {\it HNP ring} is
simply a (nonartinian) hereditary noetherian prime ring, while a {\it Dedekind
prime ring} is an HNP ring for which each nonzero ideal $I$ is {\it invertible}
in the sense that there exists a subbimodule $I^{-1}$ of the Goldie quotient
ring $Q=Q(R)$ of $R$ such that $II^{-1}=I^{-1}I=R$. Various equivalent
definitions are given in \cite{\McRob, \S5.2}.

In fact the dichotomy is stronger than the one suggested above:

\proclaim{Theorem~1} Let $R$ be a Dedekind domain that is a finitely generated 
algebra over an uncountable, algebraically closed field $k$. Then $R$ is
either simple or commutative. \endproclaim

The finite generation assumption can be further  weakened to the assertion that
$\dim_kR<\card k$ (see Theorem~6), but some such  assumption is clearly
necessary. Indeed, in the opposite direction one has the following result from
\cite{\GS}:  Suppose that $R$ is a noncommutative UFD in the sense
that $R$ is a noetherian domain for which each height one prime ideal $P$ 
satisfies (i) there exists $p\in P$ with $P=pR=Rp$  and (ii) $R/P$ is a domain. 
These exist
in profusion; for example, take the  enveloping algebra of any finite
dimensional complex nilpotent Lie algebra. The set $\C =\bigcap\{R\smallsetminus 
P\}$, where $P$ runs through the height one prime ideals of $R$,
 is localizable in $R$.  If $R$ is not 
commutative, then the localization $T=R_\C$ is necessarily a principal ideal
domain and hence a Dedekind domain \cite{\GS, Corollary~1}.
However, each height one prime ideal of $R$ induces a height one prime ideal
 of $T$.

We give a number of variants of Theorem~1. For example, an HNP domain $R$ which
is  finitely generated over an uncountable algebraically closed field $k$ 
either satisfies a polynomial identity (PI) or has no proper invertible
ideals.  To help place this in context,  recall that an HNP ring $R$ is
Dedekind prime if and only if it has no nonzero 
proper idempotent ideals \cite{\McRob,
Proposition~5.6.3}.

The key to proving all these results rests upon the following result that neither 
 requires   $R$ to be HNP nor requires  any hypotheses on $k$.

\proclaim{Theorem~2} Let $R$ be a prime noetherian algebra over a field
$k$ and $N$ a proper invertible ideal
of $R$. If $\dim_k(R/N)<\infty$, then $R$ is a PI~ring. \endproclaim

\head Proof of Theorems~1 and~2 \endhead

We first prove Theorem~2, for which we need some notation.
If $J$ is an ideal in a ring $R$, write $\C_R(J)$ for  the set of elements
of $R$ which become regular  modulo $J$. The ideal $J$ is said to be (right
and left) {\it localizable\/} if $\C_R(J)$ is a (right and left)
denominator set in $R$, in which case the localization $R[\C_R(J)^{-1}]$ is
denoted $R_J$.

\demo{Proof of Theorem~2} 
Let $J=\sqrt{N}$ denote the radical of $N$.
By \cite{\Jat, Proposition~3.3.18}, $N$ has the (right and left) Artin-Rees 
property.
 Since $J^r\subseteq N$, for some $r$, it follows that
$J$ also has the Artin-Rees property. For any $n\geq 1$, the ring
 $R/J^n$ is artinian and so $\C_{R/J^n}(J/J^n)$ is just the 
set of units in $R/J^n$. Thus  the hypotheses of \cite{\McRob,
Proposition~4.2.10} are  satisfied and, by that result, $J$ is
(right and left)  localizable.
By  \cite{\Jat, Theorem~3.2.3}, 
the localization $S= R_N$  has Jacobson radical $J(S)= JS=
SJ$. Also, $S/NS \cong R/N$. 

The 
$\C_R(J)$-torsion submodule  of $R_R$ is an ideal $I$ that contains no regular
elements. Since    $R$ is prime, this forces $I=0$.
The same result holds on the left and so
$\C_R(J)$ consists of regular elements. Thus $S$ can be identified with
a prime noetherian subring of  $Q=Q(R)$. Now
$SN^{-1}S$ is an $S$-subbimodule of $Q$ for which $(SN^{-1}S)(NS)=
S=(NS)(SN^{-1}S)$. Thus, $NS$ is an
invertible  ideal of $S$. It suffices to show that $S$ is PI
and so, replacing $R$ by $S$, we may assume that $R$ is semilocal
 with $J(R)=\sqrt{N}$.

Pick a regular element $a\in N$. 
By \cite{\Len, Theorem~2.4}, $R$   has (Rentschler-Gabriel) Krull dimension one,
and so  $R/aR$
has finite length. Hence $aR\supseteq J(R)^s\supseteq N^s$, for some $s$.  As $R/N$ 
is
finite dimensional, so are $R/N^s$ and $R/aR$. On the other hand, 
 $J(R)^t\not=J(R)^{t+1}$ for any $t$ and so 
 $\bigcap_{n\geq 0} a^nR\subseteq 
\bigcap_{n\geq 0} J(R)^n=0$. Thus \cite{\FS, Lemma~1.1}
implies that $R$ satisfies the standard identity $s_{2d}$,
 for $d=\dim_k (R/aR)$.\qed\enddemo

The following weak version of Theorem~1 is  an easy
 consequence of Theorem~2.
A field $F$ is called {\it separably closed} if it has no finite dimensional
separable extensions. The algebraic closure of a field $F$ will always be
written $\overline{F}$.
The assumptions of the next result are rather artificial
since we will subsequently need two particular  cases of it.

\proclaim{Corollary~3} Let $R$ be an  algebra over a
separably closed field $k$, such that $\dim_k R < \card k$ and
$R\otimes_k\overline{k}$ is noetherian.
Then:

{\rm (1)} If $R$ is HNP then either $R$ is a PI~ring or $R$ has no proper 
invertible ideals.

{\rm(2)} If $R$ is a Dedekind prime ring, then $R$ is either PI or simple.
\endproclaim

\demo{Proof} (1) Assume that $R$ has proper invertible ideals and let $N$ be
a maximal one. Then $N$ is semiprime \cite{\EisRo, Theorem~2.6}
and so  $R/N
\cong \bigoplus_{i=1}^t M_{n_i}(D_i)$ for some
$n_i$ and  division $k$-algebras $D_i$. 
Since $\dim_k D_i \le \dim_k R <\card k$,  \cite{\Pas,
Lemma~7.1.2} implies that each
$D_i$ is algebraic over $k$. Suppose that $D_i$ is not commutative
and let $Z_i=Z(D_i)$ denote the centre of $D_i$.
By the Koethe-Noether-Jacobson Theorem \cite{\Rwn, Exercise~2.8, p.312}
there exists $f\in D_i\smallsetminus Z_i$ such that 
$f$ is separable over $Z_i$.  
Since $k$ is separably closed, $Z_i(f)/k$ and hence $Z_i(f)/Z_i$   is
purely inseparable, giving a contradiction. 
Thus each $D_i$ is commutative.

Now $(R/N)\otimes_k\overline{k}$ is noetherian and hence so is 
 each $D_i\otimes_k\overline{k}$. Since $D_i\subseteq \overline{k}$, 
 faithful flatness implies that $D_i\otimes_kD_i$ is noetherian, 
 which by \cite{\Va, Theorem~11}
 implies that $\dim_kD_i<\infty$.
 Thus  Theorem~2 implies that $R$ is PI. 

(2)  If $R$ is not simple, then it has a proper invertible
ideal and so part (1) applies. \qed\enddemo

In order to obtain Theorem~1,  we need to understand PI
algebras satisfying the hypotheses of the corollary.

\proclaim{Lemma~4} Let $R$ be a  PI~ring over an algebraically
closed field $k$, such that $\dim_k R < \card k$. 
If $R$ is an HNP domain, then $R$ must be
commutative. \endproclaim

\demo{Proof}  By \cite{\McRob, Theorem~13.9.16 and Definition~5.3.5},  the
centre $Z=Z(R)$ is a Dedekind domain and $R$ is a finitely generated
$Z$-module.  Assume for the moment that
$Z$ has transcendence degree at most one over $k$ and let $F$ denote the field
of fractions of $Z$. 
Then the Brauer group of
$F$ is trivial by Tsen's Theorem \cite{\Pie, Corollary~a, p.376}.
Since $Q=Q(R)$  is a central simple algebra with centre $F$,
this forces $Q\cong M_n(F)$ for some $n$. But $Q$ is a
division ring, so $n=1$ and $Q=F$.
 Thus $Q$ and $R$ are commutative.  
 
 Thus, it remains to prove:\enddemo

\proclaim{Sublemma~5}
Let $Z$ be a commutative noetherian domain of Krull dimension one that is an
 algebra over an
algebraically closed field $k$ with  $\aleph=\dim_k Z < \card k=\beth$.
Then $Z$ has transcendence degree at most one over $k$.
\endproclaim

\demo{Proof}
If the sublemma fails, pick a polynomial subring $k[x,y]\subseteq Z$.
We will need the following observation.

\itemitem{($\dagger$)} Let $F[w]$ be a polynomial extension of a field $F$.
If $\lambda_i$, for $ i\in I$, are distinct elements of $F$, then the set
$\{(w-\lambda_i)^{-1} : i\in I\}\subset F(w)$ is linearly independent over $F$. 
(See, for example, the proof of \cite{\Pas, Lemma~7.1.2}.)

By ($\dagger$), the set $I=\{\lambda\in k: (x-\lambda)Z\not= Z\}$
has cardinality $\beth$. If $\lambda\in I$, then $\overline{Z}=Z/(x-\lambda)Z$ 
is a nonzero  artinian ring. By    \cite{\Pas,
Lemma~7.1.2}, again, each factor field of  $ \overline{Z}$
 is  equal to $k$ and so $\overline{Z}$
 is finite dimensional over $k$.
 This implies that the powers  of $y$ become linearly dependent in 
  $\overline{Z}$ and so there exists 
  $f_\lambda\in k[y]^*=k[y]\smallsetminus\{0\}$
  such that $f_\lambda\in (x-\lambda)Z$. In particular, 
  $z_\lambda=(x-\lambda)^{-1}f_\lambda\in Z$. 
  Let $Z(y)\subseteq Q(Z)$ denote the localization of $Z$ at
  $k[y]^*$. 
  Then    ($\dagger$) implies that  $\{z_\lambda : \lambda\in I\}$ 
  is linearly independent over $k(y)$. On the other hand, 
  $\dim_{k(y)}Z(y)\leq \aleph<\beth=\card I$, giving the required contradiction.
  This completes the proof of both Sublemma~5 and Lemma~4.
\qed\enddemo

We note that Lemma~4 also holds if the hypothesis ``$\dim_k R<\card k$''
is replaced by ``$R$ is a finitely generated $k$-algebra''.
The proof, which uses the Artin-Tate Lemma, is left to the interested reader.
 
Combining Lemma~4 with Corollary~3 gives the following generalization of
Theorem~1.

\proclaim{Theorem~6}  
Let $R$ be a domain over an
algebraically closed field $k$, such that $\dim_k R < \card k$.
Then:

{\rm (1)} If $R$ is HNP but not commutative, then $R$ has no proper 
invertible ideals.

{\rm (2)} If $R$ is a Dedekind domain, then $R$ is either commutative or
 simple. \qed
\endproclaim

\proclaim{Corollary~7} Let $R$ be a Dedekind prime ring over an
algebraically closed field $k$, such that $\dim_k R < \card k$.
Then $R$ is either simple or Morita equivalent to a commutative domain.  
\endproclaim

\demo{Proof} By \cite{\McRob, Proposition~5.2.12}, $R$ is Morita 
equivalent to a Dedekind domain. \qed\enddemo

There is  one dichotomy in the literature, due to L.~W.~Small, with a 
 flavour similar to that of Theorem~6: A prime,
noetherian, finitely generated complex algebra of Krull dimension 
one is either
primitive or satisfies a polynomial identity  \cite{\FS, p.251}.

\head Dedekind prime algebras over non-algebraically closed fields
\endhead

Theorem~1 obviously fails if we remove the hypothesis that the base field
be algebraically closed---just consider the principal ideal domain
$\HH[x]$, viewed as an algebra over $\RR$. However, this algebra is PI,
which suggests that a version of Corollary~3 might still hold. We prove
one such result in this section.

\proclaim{Lemma~8} Let $R$ be a semihereditary algebra over a field $k$.
If $F\supseteq k$ is a separable algebraic field extension, then
$R\otimes_k F$ is semihereditary. \endproclaim

\demo{Proof} Since $R\otimes_k F$ is a directed union of subalgebras
$R\otimes_k F'$ where $F'\supseteq k$ is a finite separable field
extension, we may assume that $F$ is finite over $k$. Then by \cite{\Pas,
Lemma~7.2.3}, the algebra $S= R\otimes_k F$ is {\it
relatively projective with
respect to $R$}, in the sense that
 any short exact sequence of $S$-modules which
splits as a sequence of $R$-modules also splits as a sequence of
$S$-modules.

Now consider a finitely generated right ideal $I$ of $S$ and write $I
\cong S^n/K$ for some $S$-submodule $K$ of the free $S$-module $S^n$. Since
$S_R$ is free and $I_R$ is finitely generated, $I_R$ must be projective,
and so
$K$ is an $R$-module direct summand of $S^n$. By relative projectivity,
$K$ is also an $S$-module direct summand of $S^n$, whence $I_S$ is
projective. Therefore $S$ is right semihereditary. By symmetry, it is
also left semihereditary. \qed\enddemo

\proclaim{Theorem~9} Let 
$R$ be an algebra over a field $k$ such that $R\otimes_k \kbar$ is noetherian
and $\dim_k R < \card k$. Then:

{\rm (1)} If $R$ is HNP then either
$R$ is a PI~ring or $R$ has no proper invertible ideals.

{\rm (2)} If $R$ is a Dedekind prime ring then 
$R$ is either PI or  simple.
\endproclaim

\demo{Proof} Part (2) follows immediately from part (1), so assume that $R$ is
HNP. If $\widetilde{k}$ denotes the separable closure of $k$, then faithful
flatness implies that  $\Rtilde=R\otimes_k\ktilde$ is also noetherian. By
Lemma~8, $\Rtilde$ is semihereditary and  therefore hereditary. Moreover,
$\Rtilde$ is semiprime by \cite{\Yam, Proposition~1.12}. We identify $R$ with
the $k$-subalgebra $R\otimes 1 \subseteq \Rtilde$.

By \cite{\McRob, Theorem~5.4.6}, $\Rtilde= R_0 \oplus \cdots\oplus R_n$
where $R_0$ is artinian and $R_1,\dots,R_n$ are HNP. Since $\Rtilde$ is
semiprime, $R_0$ must be semisimple. 
The regular elements of
$R$ remain regular in $\Rtilde$, and  they form a denominator set in
that ring. Hence, there are natural inclusions $\Rtilde \subseteq
Q\otimes_k \ktilde \subseteq R_0\oplus Q_1\oplus \cdots\oplus Q_n$ where
$Q=Q(R)$  and $Q_\ell=Q(R_\ell)$ for $\ell\ge1$. Note that
$$\dim_{\tilde{k}} R_\ell  \le \dim_k R <\card k\le \card \ktilde$$ and so
Corollary~3 applies to $R_\ell$ for $\ell\ge 1$. 

Suppose that $R$ has a proper invertible ideal $I$. Then $I$ induces an ideal
$J= \bigoplus_{\ell=0}^n J_\ell\subseteq \Rtilde$, and $I^{-1}$ induces an
$\Rtilde$-subbimodule $J'\subseteq Q\otimes_k \ktilde$ such that $JJ'=J'J=
\Rtilde$. It follows that each $J_\ell$ is an invertible ideal in $R_\ell$.
Moreover, using \cite{\McRob, (5.6.2)}, $\bigcap_{j\geq 0}J_\ell^j\subseteq
\bigcap_{j\geq 0}(I\otimes\widetilde{k})^j=0$ for each $\ell$. If $\ell=0$, the
first of these observations implies that $J_0=R_0$, and the second then implies
that $R_0=0$.

On the other hand,  $J_\ell$ is a proper  invertible ideal in $R_\ell$ for
$\ell\geq 1$. Clearly   $R_\ell\otimes_{\tilde{k}}\overline{k}$ is a
summand of
$\widetilde{R}\otimes_{\tilde{k}}\overline{k}=R\otimes_k\overline{k}$ and
so  $R_\ell\otimes_{\tilde{k}}\overline{k}$ is noetherian. Thus the
hypotheses of  Corollary~3 are satisfied and each  $R_\ell$  is~PI. Therefore
$R$ is~PI. \qed\enddemo

As we mentioned in the introduction, the examples from \cite{\GS} show that 
the theorem  fails badly if $\dim_k R = \card k$. Similarly it can fail
badly if $R\otimes_k \kbar$ is not noetherian. Here is a typical example.
Let $D$ denote the Krull division ring; thus $D=\bigotimes_{\Bbb Q} D_p$, where
for each prime $p\in \Bbb Z$, $D_p$ is a $p^2$-dimensional, $\Bbb Q$-central
division ring (see, for example, \cite{\ssw, p.221}). Let 
$k=\Bbb Q(x_i : i\in I)$ be a purely transcendental field extension of 
$\Bbb Q$ of cardinality  $\geq \aleph_1$.
 Then, $E=D\otimes_{\Bbb Q}k$ is a division ring and  $R=E[y]$ is
a Dedekind domain for which the conclusion of Theorem~9 fails. 
Indeed, for each prime $p$ there exists  $f_p\in \Bbb Q[y]$ such that
$R/f_pR\cong M_{p}(F_p)$, for the appropriate division ring
 $F_p$ algebraic over $k$. Thus, $R$ is not PI
but, of course, each $f_pR$ is an invertible ideal and  $\dim_kR=\aleph_0<\card
k$. Examples of other HNP rings with distinctive  ideal structures can be found
in \cite{\ssw} and \cite{\sww}. To give one example, let 
  $H$ be the ring  constructed in
\cite{\sww, Theorem~2.2}. Then $S=H\otimes_{\Bbb Q}k$ will be a 
subring of $R$ such that 
 $S$ is an HNP ring with infinitely many idempotent ideals.

The first example from the last paragraph also shows that Theorem~2 will
fail if we weaken the assumption that $\dim_k (R/N)<\infty$ to  ``$\dim_k
(R/N)<\card k$'' or even to  ``$R/N$ is algebraic over the uncountable field
$k$.''

\enddocument